\documentclass[12pt]{amsart}
\def\ds{\displaystyle}
\newtheorem{theorem}{Theorem}
\newtheorem{proposition}[theorem]{Proposition}
\newtheorem{lemma}[theorem]{Lemma}
\def\CC{{\mathbb C}}
\def\DD{{\mathbb D}}

\def\Re{\operatorname{Re}}
\def\Im{\operatorname{Im}}
\def\Id{\operatorname{Id}}
\def\co{\operatorname{co}}
\def\dist{\operatorname{dist}}

\title{Entire curves avoiding given sets in $\CC^n$}
\author{Nikolai Nikolov and Peter Pflug }
\address
{Institute of Mathematics and Informatics\\ Bulgarian Academy of
Sciences\\ 1113 Sofia, Bulgaria} \email{nik@math.bas.bg}
\address{Carl von Ossietzky Universit\"at Oldenburg\\
Fakult\"at V, Institut f\"ur Mathematik\\ Postfach 2503\\ D-26111 Oldenburg,
Germany} \email{pflug@mathematik.uni-oldenburg.de}

\begin{document}


\begin{abstract}
Let $F\subset\CC^n$ be a proper closed subset of $\CC^n$ and
$A\subset\CC^n\setminus F$ at most countable ($n\ge 2$). We give conditions on $F$ and $A$,
under which there exists a holomorphic immersion (or a proper holomorphic embedding)
$\varphi:\CC\to\CC^n$ with $A\subset\varphi(\CC)\subset\CC^n\setminus F$.
\end{abstract}

\maketitle

Let $F\subset\CC^n$ be a proper closed subset of $\CC^n$ and
$A\subset\CC^n\setminus F$ at most countable ($n\ge 2$). The aim of this
note is to discuss conditions for $F$ and $A$, under which there exists a
holomorphic immersion (or a proper holomorphic embedding)
$\varphi:\CC\to\CC^n$ with $A\subset\varphi(\CC)\subset\CC^n\setminus F$.
Our main tool for constructing such mappings is Arakelian's approximation
theorem (cf. \cite{GP,RR3}).

The first result is a generalization of the main part of Theorem 1
in \cite{JN}. More precisely, we prove the following result.

\begin{proposition}\label{prop1} Let $F$ be a proper convex closed set in $\CC^n,\ n\ge 2$.
Then the following statements are equivalent:

(i) either $F$ is a complex hyperplane or it does not contain any
complex hyperplane;

(ii) for any integer $k\ge 1$ and any two sets
$\{\alpha_1,\alpha_2,\dots,\alpha_k\}\subset\CC$ and
$\{a_1,a_2,\ldots,a_k\}\subset\CC^n\setminus F$, there exists a proper
holomorphic embedding $\varphi:\CC\to\CC^n$ such that
$\varphi(\alpha_j)=a_j,\ 1\le j\le k$, and
$\varphi(\CC)\subset\CC^n\setminus F$.

(iii) the same as (ii) but for $k=2$.

\end{proposition}

The equivalence of  ($i$) and ($iii$) follows from the proof of
Theorem 1 in \cite{JN}. For the convenience of the reader we repeat here the
main idea of the proof of $(iii)\Longrightarrow (i)$. Observe that condition
($iii$) implies that the Lempert function of the domain $D:=\CC^2\setminus
F$ is identically zero, i.e.
\begin{multline}{\notag}
$$
\widetilde k_D(z,w):=\inf\{\alpha\ge 0: \\ \exists f\in\mathcal O(\Delta,D):
f(0)=z,\;f(\alpha)=w\}=0,\quad z,w\in D,
$$
\end{multline}
where $\Delta$ denotes the open unit disc in $\CC$. In the case when
condition ($i$) is not satisfied we may assume (after a biholomorphic
mapping) that $F=A\times\CC$, where the closed convex set $A$, properly contained in $\CC$,
contains at least two points. Applying standard properties of $\widetilde
k$, we have $\widetilde k_D(z,w)=\widetilde
k_{\CC\setminus A}(z',w')$, where $(z,w)=((z',z''),(w',w''))\in D$. Since
$\widetilde k_{\CC\setminus A}$ is not identically zero we end with a
contradiction.

Hence, we only have to prove the implication
$(i)\Longrightarrow(ii)$.

\begin{proof} For simplicity of notations we shall consider only the case $n=2$.

If $F$ is a complex line, we may assume that $F=\{z_2=0\}$. Considering an
automorphism of the form $(z_1,z_2)\to(z_1e^{\gamma z_1z_2},z_2e^{-\gamma
z_1z_2})$ for a suitable constant $\gamma$, we may also assume that the
second coordinates of the given points are pairwise different. Then there
exist two one variable polynomials $P$ and $Q$ such that the mapping
$t\to(t+P(e^{Q(t)}),e^{Q(t)})$ has the required property.

Assume now that $F$ does not contain any complex line. The idea below
comes from that of Theorem 8.5 in \cite {RR1}.

First, we shall prove by induction that for any $j\le k$ there is
an automorphism $\Phi_j$ such that the set $\co(\Phi_j(F))$ does
not contain any complex line and it does not have a common point
with the set
$$
\co(G_j)\cup\{\Phi_j(a_{j+1}),\ldots,\Phi_j(a_k)\},
$$
where $G_j:=\{\Phi_j(a_1),\ldots,\Phi_j(a_j)\}$ ($\co(M)$ denotes
the convex hull of a closed set $M$ in $\CC^n$). Doing the
induction step, we may assume that $\Phi_j=\Id$. Then, since $F$
is convex and does not contain any complex line, after an affine
change of coordinates one has that (cf. \cite{D,JN})
$$
F\subset H:=\{\Re(z_1)\le -1,\Re(z_2)\le -1\},
$$
$$
\co(G_j)\subset\{\Re(z_1)\ge 0\},\ a_{j+1}\in\{\Re(z_2)\ge 0\}.
$$
In addition, we may assume that the set $A:=\{a_1,\ldots,a_k\}$ of
the given points and the strip $\{-1<\Re(z_2)<0\}$ do not have a
common point. By Arakelian's theorem (cf. \cite{GP}), for
$\varepsilon:=\min\{1,\dist(F,A)\}$ we may find an entire function
$f$ such that
$$
|f(t)-a_{j+1,1})|<\frac{\varepsilon}{2}\hbox{ if
}\Re(t)\le -1,\ |f(t)|<\frac{\varepsilon}{2}\hbox{ if }\Re(t)\ge
0
$$
and, in addition, $f(a_{j+1,2})=0$ (here, $a_{j+1,k}$ denotes the
$k$-th coordinate of the point $a_{j+1}$). Then it is easy to see
that the automorphism $\Phi_{j+1}(z_1,z_2):=(z_1+f(z_2),z_2)$ has
the required properties.

So, let $F$ be a convex set, which does not contain any complex
line and $F\cap\co(A)=\emptyset$. Then we may assume that (cf.
\cite{D,JN}) $F\subset H$, $A\subset\{\Re(z_1)\ge 1,\Re(z_2)\ge
0\}$, and, in addition, that $\Re(\alpha_j)\ge 1,\ 1\le j\le k$.

Note that there exists an entire function $g$ such that $|g(t)|\le
1$ if $\Re(t)\le -1$ and $g(a_{j,2})=\alpha_j-a_{j,1}$ (cf.
\cite{GP,RR2}; this can be proved also directly, applying a
standard interpolation process and Arakelian's theorem many
times). Then, applying the automorphism
$(z_1,z_2)\to(z_1+g(z_2),z_2)$, we may assume that
$a_{j,1}=\alpha_j$ and $F\subset\{\Re(z_1)\le 0,\Re(z_2)\le -1\}$.
Finally, we find, as above, an entire function $h$ such that
$|h(t)|<1$ on the set $\Re(t)\le 0$ and $h(\alpha_j)=a_{j,2}$.
Hence, the mapping $t\to(t,h(t))$ has the required properties (in
the new coordinates).
\end{proof}

The end of the proof shows that we may also prescribe values of finitely
many derivatives of $\varphi$ at the points of the given planar set.
\smallskip

\noindent
{\bf Open problem.} Is it true for an $F$ as in  ($ii$) of Proposition
1, that for any discrete set of points in $\CC^n\setminus F$ there exists a
proper holomorphic embedding of $\CC$ in $\CC^n$ avoiding $F$ and passing
through any of these points?
\smallskip

It is known that for any discrete set of points in $\CC^n$ there
exists a proper holomorphic embedding of $\CC$ in $\CC^n$ passing
through any of the points of this set (Proposition 2 in
\cite{FGR}; cf. also Theorem 1 in \cite {RR2} for $n\ge 3$). We
have not been able to modify the proofs of \cite{FGR,RR2} to get a
positive answer for the above question in the general case.
Nevertheless, the following result gives a positive answer to the
open problem in the case when $F$ is a complex hyperplane.

\begin{proposition}\label{prop2} If $F$ is a union of at most $n-1$ $\CC$-linearly
independent complex hyperplanes in $\CC^n$, then for any discrete set of
points in $\CC^n\setminus F$ there exists a proper holomorphic embedding of
of $\CC$ into $\CC^n$ avoiding $F$ and passing through any of these points.
\end{proposition}

The proof of Proposition \ref{prop2} will be a modification of the
one in the case when $F$ is the empty set ( see Proposition 2 in
\cite{FGR}).

The key point is the following

\begin{lemma}\label{prop3} Let $K$ be a polynomially convex compact set in $\CC^n$,
$A$ a set of finitely many points in $K$, and $H$ a union of at
most $n-1$ linearly independent complex hyperplanes in $\CC^n$.
For every $p,q\in\CC^n\setminus(K\cup H)$ and every
$\varepsilon>0$, there exists an automorphism $\varphi$ of $\CC^n$
such that $\varphi(z)=z$, $z\in H\cap A$, $\varphi(p)=q$, and
$|\varphi(z)-z|\le\varepsilon$, $z\in K$.
\end{lemma}

In view of Lemma \ref{prop3}, Proposition \ref{prop2} follows by
repeating step by step the proof of Proposition 2 in \cite{FGR}.
Starting with an embedding $\alpha_0$ whose graph avoids $H$, the
desired embedding $\alpha$ is constructed as the limit of a
sequence of embeddings $\alpha_j$ with
$\alpha_j=\varphi_j\circ\alpha_{j-1}$ $(j\ge 1)$, where the
$\varphi_j$ are automorphisms chosen by Lemma \ref{prop3}. Note
that the graph of $\alpha$ avoids $H$ by the Hurwitz theorem.
\smallskip

\noindent {\it Proof of Lemma 3.} After a linear change of
coordinates, we may assume that $H\subset\{z_1\cdots z_n=0\}$ and
that all the coordinates of the points in $B:=A\cup\{q\}\setminus
H$ are non-zero. Applying an overshear of the form
$$
w_1=z_1\exp(f(z_2,\dots,z_n)), w_2=z_2,\dots, w_n=z_n,
$$
where
$$
f(z_2,\dots,z_n):=z_2\dots
z_n(\varepsilon+\sum_{j=2}^{n}\varepsilon^jz_j)
$$
and $\varepsilon$ is small
enough, provides pairwise different products of the first $n-1$ coordinates
of the points in $B$. Repeating this argument, we may assume the same for
every $n-1$ coordinates.
\smallskip

Now, we need the following variation of Theorem 2.1 in \cite{FR}.

\begin{lemma}\label{prop4} Let $H$ be the union of at most $n-1$ linearly
independent complex hyperplanes in $\CC^n$, $D$ an open set in $\CC^n$, and
$K\subset D$ a compact set. Let $\Phi_t:D\to\CC^n$, $t\in[0,1]$, be a
$C^2$-smooth isotopy of biholomorphic maps which fix pointwise $D\cap H$
such that $\Phi_t(D\cap H)=\Phi_t(D)\cap H$. Suppose that $\Phi_0$ is the
identity map and the set $\Phi_t(K)$ is polynomially convex for every
$t\in[0,1]$.

Then $\Phi_1$ can be approximated, uniformly on $K$, by
automorpisms of $\CC^n$, which fix pointwise $H$.
\end{lemma}

For a moment, we may assume that Lemma \ref{prop4} is true. Let
$\gamma:[0,1]\to\CC^n\setminus(K\cup H)$ be a $C^2$-smooth path,
$\gamma(0)=p$, $\gamma(1)=q$. Then we apply Lemma \ref{prop4} to
the following situation:

take $\Phi_t(z)$ to be $z$ near $K$ and to be $z+\gamma(t)-p$ near
$p$, and choose a sufficiently small neighborhood $D$ of the
polynomially convex set $K\cup\{p\}$. For a sufficiently small
$\varepsilon>0$, denote by $\psi$ the corresponding automorphism
and set $\tilde r:=\psi(r)$ for $r\in B$. Let $f_1$ be the
Lagrange interpolation polynomial with
$$
f_1(\tilde r_2\dots\tilde r_n)=\frac{1}{\tilde r_2\dots\tilde
r_n} \log\frac{r_1}{\tilde r_1}
$$
for every $r\in B$. Note that the overshear
$$
\psi_1(z):=(z_1\exp(z_2\dots z_nf_1(z_2\dots z_n)),z_2,\dots,z_n)
$$
sends $\tilde r$ to the point $(r_1,\tilde r_2,\dots,\tilde r_n)$.
It is left to define in a similar way $\psi_2,\dots,\psi_n$ and to
consider the composition $\psi_n\circ\dots\circ\psi_1\circ\psi$.\
This completes the proof of Lemma \ref{prop3}.\qed
\smallskip

\noindent {\it Proof of Lemma 4.} Note that under the assumptions
of Lemma \ref{prop4}, there exists a neighborhood $U\subset D$ of
$K$ such that $U_t:=\Phi_t(U)$ is Runge for each $t\in[0,1]$
(Lemma 2.2 in \cite{FR}). We shall follow the proofs of Theorem
1.1 in \cite{FR} and Theorem 2.5 in \cite{Var}. Consider the
vector field $\ds X_t:=\frac{d}{dt}\Phi_t\circ\Phi_t^{-1}$ defined
on $U_t$. For a sufficiently large positive integer $N$ and $0\le
j\le N-1$ set
$$
X_{j,t}:=\left\{\begin{array}{ll} 0,&t\not\in[j/N,(j+1)/N]\\
X_{j/N},&t\in[j/N,(j+1)/N].\\
\end{array}\right.
$$
Note that $X_{j/N}$ vanishes on $U_{j/N}\cap H$. It is easy to see
that it can be approximated by holomorphic vector fields on
$\CC^n$ which vanish on $H$, since $U_{j/N}$ is Runge (here and
below, the approximations are locally uniformly). On the other
hand, these vector fields can be approximated by Lie combinations
of complete vector fields vanishing on $H$ (Proposition 5.13 in
\cite{Var}). Thus we may assume that $X_{j/N}$ is a Lie
combination of complete vector fields vanishing on $H$. Note that
the local flow of $\ds\sum_{j=0}^{N-1}X_{j,t}$ at time $1$ is
$h_{N-1}\circ\dots\circ h_0$, where $h_j$ is the local flow of
$X_{j/N}$ at time $\ds\frac{1}{N}$. If $N\to\infty$, then this
composition converges to the time one map $\Phi_1$ of the flow of
$X_t$. To finish the proof of Lemma \ref{prop4}, it is enough to
note that every $h_j$ can be approximated by finite compositions
of automorphisms of $\CC^n$ which fix $H$ (cf. the proof of
Theorem 2.5 in \cite{Var}).\qed
\smallskip

In this way Proposition \ref{prop2} is completely proved.
\smallskip

{\it Remark.} It is an open question whether every holomorphic vector field
in $\CC^n$, which vanishes on the set $L:=\{z_1\cdots z_n\}$, can be locally
uniformly approximated  by Lie combinations of complete vector fields
vanishing on $L$ \cite{Var}. If this would be so, then the above proof shows
that Proposition \ref{prop2} is also true for every union of linearly
independent complex hyperplanes in $\CC^n$, $n\ge 3$. To see this, choose,
for example, the starting embedding
$$
\alpha_0(\eta):=(\exp(-\eta^2),\exp(-\eta\sqrt2),\exp(\eta),\dots,\exp(\eta)).
$$ It remains an unsolved problem (for us) if there exists a proper
holomorphic embedding of $\CC$ in $\CC^2$ whose graph avoids both coordinate
axes.
\smallskip

We are also able to answer the open problem, posed after
Proposition \ref{prop1}, in the bounded case.

\begin{proposition}\label{prop5} If $K$ is a polynomially convex compact set in $\CC^n$,
then for any discrete set $C$ of points in $\CC^n\setminus K$
there exists a proper holomorphic embedding $H$ of $\CC$ in
$\CC^n$ avoiding $K$ and passing through any of these points. In
addition, for a given point $c\in C$ and $X\in\CC^n\setminus\{0\}$
we can choose $H$ such that $H'(H^{-1}(c))=X$. In particular, the
Lempert function and the Kobayashi pseudometric of $\CC^n\setminus
K$ vanish.
\end{proposition}

\begin{proof} The proof is a modification of the one of Proposition 2 in \cite{FGR}.

We may assume that $X=(1,0,\dots,0)$ and that $K$ does not
intersect the first coordinate axis. Note that there exists a
smooth non-negative plurisubharmonic exhaustion function $\varphi$
on $\CC^n$ that is strongly plurisubharmonic on $\CC^n\setminus K$
and vanishes precisely on $K$ (cf. \cite{Ca}). For any $\epsilon
>0$, put
$$
G_\varepsilon:=\{z\in\CC^n: \phi(z)<\varepsilon\} \text{ and }
K_\varepsilon:=\{z\in\CC^n: \phi(z)\le\varepsilon\}. $$ In
particular, $K_\varepsilon$ is polynomially convex. By Sard's
theorem we may choose a strictly decreasing sequence
$(\varepsilon_j)_{j\ge 0}$, bounded from below by a positive
constant, such that the boundary of $G_j:=G_{\varepsilon_j}$ is
smooth for any $j$ and $K_0:=K_{\varepsilon_0}$ does not intersect
the first coordinate axis. In particular, $K_j:=K_{\varepsilon_j}$
has finitely many connected components.

\noindent {\it Claim.} $K_j\subset\psi_j(K_{j-1})$ for any
automorphism $\psi_j$ of $\CC^n$  which is closed enough to the
identity map on $K_{j-1}$.

Let now $C=(\alpha_l)_{l\ge 1}$ with $\alpha_1=c$. Set
$H_0(\zeta)=(\zeta,0,\dots,0)$ and $\rho_0=0$. In view of the
claim and the proof of Proposition 2 in \cite{FGR}, for any $j\ge
1$ we may find by induction numbers $\rho_j\ge\rho_{j-1}+1$,
$\zeta_j\in\CC$, and an automorphism $\psi_j$ such that for
$H_j=\psi_j\circ H_{j-1}$ one has:

(a) $H'_j(\zeta_1)=X$ and $H_j(\zeta_l)=\alpha_l, 1\le l\le j$;

(b) $|H_j(\zeta)|>|\alpha_j|-1$ if $|\zeta|\ge\rho_j$ and $\ds
K_j\subset\{z\in\CC^n:|z|\le|\alpha_j|-\frac{1}{2}\}$;

(c) $|H_j(\zeta)-H_{j-1}(\zeta)|\le\delta_j\le2^{-j}$ if
$|\zeta|\le\rho_j$;

(d) $H_j(\CC)\cap K_j=\emptyset$.

It is easy to check that the limit map $\ds
H:=\lim_{j\to\infty}H_j$ exists and that it has the required
properties except properness. The last one can be provided by the
choice of $\delta_j$. Note that the only modifications that have
to be made in the proof of Proposition 2 in \cite{FGR} are the
choice of the $\psi_j$ with the additional property
$\psi'_j(\zeta_1)$ to be the unitary matrix and the replacing of
the set $\ds F:=\{z\in\CC^n:|z|\le|\alpha_j|-\frac{1}{2}\}\cup
H_{j-1}\{|\zeta|\le\rho\}$ by the set $F:=K_j\cup
H_{j-1}\{|\zeta|\le\rho\}$ if $\ds
K_j\not\subset\{z\in\CC^n:|z|\le|\alpha_j|-\frac{1}{2}\}$.

\noindent {\it Proof of the claim.} Since $K_j$ has finitely many connected
components $K_{j,1},\dots,K_{j,m}$, we have that $\dist(K_j,\partial
K_{j-1})>0$. Then we find an $r>0$ with $\dist(K_j,\partial K_{j-1})>r$ and
some ball $B_l$ with radius $r$ belonging to $K_{j,l}$, $1\le l\le m$. It
follows that $K_j\subset\psi_j(K_{j-1})$, if $$ \max\{|\psi_j(z)-z|:z\in
K_{j-1}\}\le r. $$ Indeed, suppose the contrary, i.e., $\psi_j(a)\in K_j$
for some $a\not\in K_{j-1}$. We may assume that $\psi_j(a)\in K_{j,1}$.
Denote by $b_1$ the image of the center of $B_1$ under $\psi_j$. Then there
exists a path $\gamma$ in $K_{j,1}$ joining $\psi_j(a)$ and $b_1$. Note that
$\psi_j^{-1}(\gamma)\cap\partial K_{j-1}\neq\emptyset$. If
$c\in\psi_j^{-1}(\gamma)\cap\partial K_{j-1}$, then $\psi(c)\in K_j$. Hence
$r\ge|\psi_j(c)-c|\ge\dist(K_j,\partial K_{j-1})$; a contradiction.
\end{proof}

Note that if $F$ is a proper subset in $\CC^2$ such that for any point in
$\CC^2\setminus F$ there exists a non-constant entire curve
$\gamma:\CC\to\CC^2\setminus F$ which passes through this point, then the
interior of $F$ is pseudoconvex, since
$\CC^2\setminus\overline{\gamma(\CC)}$ is pseudoconvex \cite{U}. Moreover,
if $F$ is compact and for any point $a\in\CC^2\setminus F$ there exists a
proper holomorphic mapping $\varphi:\CC\to\CC^2$ with
$a\in\varphi(\CC)\subset\CC^2\setminus F$, then $F$ is rational convex
\cite{GR}. The same does not holds in higher dimensions. For example, if $F$
and $G$ are two proper closed subsets of $\CC^k$ and $\CC^l$, respectively,
then for any point in $\CC^{k+l}\setminus(F\times G)$ there exists a proper
holomorphic embedding of $\CC$ in $\CC^{k+l}$ avoiding $F\times G$ and
passing through this point.

The next proposition is in the spirit of the above remark and it generalizes
Proposition 1 in \cite{Nik}.

\begin{proposition}\label{prop6} If $F$ and $G$ are two sets in $\CC^k$ and
$\CC^l$,
respectively, then for any countably set $C$ of points in
$\CC^{k+l}$ with $\dist(C,F\times G)>0$ there exists a holomorphic
immersion of $\CC$ in $\CC^{k+l}$ avoiding $F\times G$ and passing
through any point of $C$.
\end{proposition}

\begin{proof} The idea for the proof comes from the one of Theorem 2 in
\cite{RR2}. For any point $c$ in $\CC^{k+l}$ denote by $c'$ and $c''$ its
projections onto $\CC^k$ and $\CC^l$, respectively. Set
$\varepsilon:=\dist(C,F\times G)>0$, $C':=\{c\in C:\dist(c',\CC^k\setminus
F)\ge\varepsilon\}$ and $C'':=C\setminus C'$. We may assume that both sets
are infinite and enumerate them, i.e. $C'=(a_j)_{j\ge 0}$ and
$C''=(b^j)_{j\ge 0}$. Denote by $\DD_n(c,r)$ the polydisc in $\CC^n$ with
center at $c$ and radius $r$. Note that
$\DD_k(a'_j,\varepsilon)\subset\CC^k\setminus F$ and
$\DD_l(b''_j,\varepsilon)\subset\CC^l\setminus G$ for any $j\ge 0$. Define
$$ A_j:=\{z\in\CC:\Re(z)\le -3,|\Im(z)-7j|\le 3\},\ j\ge 1, $$ $$
A_0:=\{z\in\CC:\Re(z)\ge -1\}\setminus
\cup_{j=1}^\infty\{z\in\CC:\Re(z)>5,|\Im(z)-7j|<1\}. $$ Choose a number
$t\in(0,1)$ such that $$ t\exp({\root3\of
x}-{\root4\of{x+2}})\ge4e(1-t){\root3\of{(x+2)^4}},\ x\ge 0. $$ For $1\le
m\le k$, combining the extensions of Arakelian's theorem in  \cite{GP,RR2}
gives an entire function $f_m$ such that $$
|f_m(z)-a_{0,m}-\frac{\varepsilon}{2}t\exp(-{\root4\of{z+2}})|<
\frac{\varepsilon}{2}(1-t)\exp(-{\root3\of{|z|}}),z\in A_0, $$ $$
|f_m(z)-a_{j,m}-\frac{\varepsilon}{2}t\exp(-{\root4\of{-z-2}})|<
\frac{\varepsilon}{2}(1-t)\exp(-{\root3\of{|z|}}),z\in A_j,j\ge 1, $$ $$
f_m(2)=a_{0,m},\;f_m(-2)=b_{0,m},\;f_m(-7+i7j)=a_{j,m},\;f_m(7+i7j)=b_{j,m}
$$ for $j\ge 1$ ($\root4\of z$ is the branch with $\root4\of 1=1$ and
$c_{j,m}$ denotes the m-th coordinate of the point $c_j$). Note that
$|f_m(z)-a_{j,m}|<\varepsilon$ if $z\in A_j$. For $k+1\le m\le k+l$ we
choose analogously an entire function $f_m$ such that $$
|f_m(z)-b_{0,m}-\frac{\varepsilon}{2}t\exp(-{\root4\of{-z-2}})|<
\frac{\varepsilon}{2}(1-t)\exp(-{\root3\of{|z|}}),-z\in A_0, $$ $$
|f_m(z)-b_{j,m}-\frac{\varepsilon}{2}t\exp(-{\root4\of{z+2}})|<
\frac{\varepsilon}{2}(1-t)\exp(-{\root3\of{|z|}}),-z\in A_j,j\ge 1, $$
$$f_m(2)=a_{0,m},\;f_m(-2)=b_{0,m},\;f_m(-7+i7j)=a_{j,m},\;f_m(7+i7j)=b_{j,m}
$$ for $j\ge 1$. Then the mapping $(f_1,\dots,f_{k+l})$ will have the
required properties if it is non-singular. To see this, note that applying
the triangle inequality and the Cauchy inequality gives $$
\Bigl|\frac{\varepsilon t}{8\root3\of{|z+2|^4}}\exp(-\root4\of{|z+2|})\Bigr|
-|f'_m(z)|<\frac{\varepsilon}{2}(1-t)\exp(1-{\root3\of{|z|}}) $$ for $1\le
m\le k$ and $$ z\in E_0:=\{z\in\CC:\Re(z)\ge 0\}\setminus
\bigcup_{j=1}^\infty\{z\in\CC:\Re(z)>4,|\Im(z)-7j|<2\}. $$ Then the choice
of $t$ shows that $f'_m(z)\neq 0$ if $1\le m\le k$ and $z\in E_0$; a similar
argument gives that $f'_m(z)\neq 0$ if $$ z\in
E:=\bigcup_{j=1}^\infty\{z\in\CC:\Re(z)\le -4,|\Im(z)-7j|\le 2\}. $$ We
obtain analogously that $f'_m(z)\neq 0$ if $k+1\le m\le k+l$ and $-z\in
E_0\cup E$, which implies that the mapping is non-singular.
\end{proof}

Note that, in general, the mapping in Proposition \ref{prop6} cannot be chosen to be proper.
For example, let $F:=\CC\setminus\DD_1(0,1)$ and let
$f:=(f_1,f_2)$ be a proper holomorphic map of $\CC$ in $\CC^2$
which avoids $F\times F$. Choose an $R$ such that
$\max\{|f_1(z)|,|f_2(z)|\}\ge 2$ for $|z|>R$. Assume that $f_1$ is
not a polynomial. Then by Picard's theorem there is a point
$a\in\CC$, $|a|>R$, with $|f_1(a)|=1$. Thus $|f_2(a)|\ge 2$. On
the other side, using that $f(\CC)\cap(F\times F)=\emptyset$
implies that $|f_1(a)|<1$, a contradiction. In conclusion, one of
the functions $f_1$ and $f_2$ is a polynomial and the other one is
a constant smaller than 1.

It follows from Proposition \ref{prop6} that if $F$ and $G$ are two closed
proper subsets of $\CC^k$ and $\CC^l$, respectively, then the
Lempert function of $\CC^{k+l}\setminus(F\times G)$ vanishes. The
next proposition implies that the same holds for the Kobayashi
pseudometric.

\begin{proposition} \label{prop7} If $F$ and $G$ are two proper closed sets in $\CC^k$ and
$\CC^l$, respectively, then for any point
$c\in\CC^{k+l}\setminus(F\times G)$ and any vector $X\in\CC^{k+l}$
there exists a holomorphic mapping of $\CC$ in
$\CC^{k+l}\setminus(F\times G)$ with $f(0)=c$ and $f'(0)=X$.
\end{proposition}

\begin{proof} We may assume that $c'\in\CC^k\setminus F$ and
$\DD_l(0,1)\subset\CC^l\setminus G$. The statement is trivial if
$X'=0$. Otherwise, we may assume $c'=0$ and the ball in $\CC^k$
with center at the origin and radius $(e+1)\sqrt k$ belongs to
$\CC^k\setminus F$. After a unitary transformation of $\CC^k$ we
may also assume that $X'=(r,\dots,r)$ for some $r>0$. Note that
$\DD_k(0,e+1)\subset\CC^k\setminus F$ and if $|e^{rz}-1|\ge e+1$,
then $\ds\Re(z)\ge\frac{1}{r}$. By Arakelian's theorem, there
exists an entire function $f_m$ such that $f_m(0)=0,\ f'_m(0)=X_m$,
and $|f_m(z)|<1$ if $\ds\Re(z)\ge\frac{1}{r},\ k+1\le m\le k+l$.
Setting $f_m(z):=e^{rz}-1$ for $1\le m\le k$ implies that the
mapping $(f_1,\dots,f_{k+l})$ has the required properties.
\end{proof}

\noindent {\bf Acknowledgments.} A part of this note was prepared during the
stay of the first author at the University of Oldenburg (January, 2003),
supported by DFG, and during the stays of both authors at the Jagiellonian
University in Krak\'ow (February, 2003) with support for the second author
by DFG. We like to thank all these institutions.

Proposition \ref{prop2} was obtained after the first version of this note
was accepted for publication. We would like to thank our colleagues Franc
Forstneri\v c and Josip Globevnik for extremely helpful discussions about
this proposition.

\end{document}